\documentclass[10pt]{amsart}
\usepackage{latexsym,amssymb} 
\usepackage{xy}
\xyoption{all}

\theoremstyle{plain}
  \newtheorem{theorem}{Theorem}[section]
  \newtheorem{proposition}[theorem]{Proposition}
  \newtheorem{lemma}[theorem]{Lemma}
  
  \newtheorem{conjecture}[theorem]{Conjecture}
\theoremstyle{definition}
  \newtheorem{definition}[theorem]{Definition}
  \newtheorem{example}[theorem]{Example}
  
  \newtheorem{question}[theorem]{Question}
 \theoremstyle{remark}
  \newtheorem{remark}[theorem]{Remark}

\numberwithin{equation}{section}

\def\ZZ{{\mathbb Z}}
\def\QQ{{\mathbb Q}}

\def\CC{{\mathbb C}}
\def\HH{{\mathbb H}}

\def\ff{{\mathbf f}}
\def\thetas{{\mathbf{\Theta}}}

\def\symm{{\mathfrak S}}
\def\antichains{{\mathfrak{An}}}

\def\im{\mathrm{im}}
\def\re{\mathrm{re}}

\def\Cat{\mathrm{Cat}}
\def\Nar{\mathrm{Nar}}
\def\Cent{\mathrm{Cent}}

\def\Sym{\mathrm{Sym}}
\def\Hilb{\mathrm{Hilb}}
\def\Gal{\mathrm{Gal}}
\def\triv{\mathrm{triv}}

\newcommand{\qbin}[2]{\left[ \begin{matrix} #1 \\ #2 \end{matrix} \right]_q}

\begin{document}

\title[Cyclic sieving of noncrossing partitions]
{Cyclic sieving of noncrossing partitions for complex reflection groups}

\author{David Bessis}
\address{DMA - \'Ecole normale sup\'erieure\\
45 rue d'Ulm\\
75230 Paris cedex 05\\
France}
\email{david.bessis@ens.fr}

\author{Victor Reiner}
\address{School of Mathematics\\
University of Minnesota\\
Minneapolis, MN 55455, USA}
\email{reiner@math.umn.edu}

\thanks{Second author supported by NSF grant DMS--9877047.}



\begin{abstract}
We prove an instance of the cyclic sieving phenomenon,
occurring in the context of noncrossing parititions for well-generated
complex reflection groups.
\end{abstract}

\maketitle

\section{Introduction}

  Our goal is Theorem~\ref{main-result} below,
whose terminology is explained briefly here, and more fully in the next section.

Let $V=\CC^n$ and let $W \subset GL(V)$ be a finite irreducible complex reflection group, that is,
$W$ is generated by its set $R$ of reflections.  Shephard and Todd \cite{ShephardTodd} classified
such groups, and used this to show that when $W$ acts on the symmetric algebra $S:=\Sym(V^*)$, its
invariant subring $S^W$ is itself a polynomial algebra, generated
by homogeneous polynomials $f_1,\ldots,f_n$ whose 
degrees $d_1 \leq d_2 \leq \cdots \leq d_n$ are uniquely determined;  this was
also proven without using the classification by Chevalley \cite{Chevalley}.  We will
assume further that $W$ is {\it well-generated} in the sense that it 
can be generated by $n$ of its reflections.
Define the Coxeter number $h:=d_n$ and the $W$-$q$-{\it Catalan number}
\begin{equation}
\label{W-q-Catalan-definition}
\Cat(W,q):=\prod_{i=1}^n \frac{[h+d_i]_q}{[d_i]_q}
\end{equation}
where $[n]_q:=1+q+q^2+\cdots+q^{n-1}=\frac{q^n-1}{q-1}.$
Let $c$ be a regular element of $W$ in the sense of Springer \cite{Springer},
of order $h$;  such a regular element of order $h$ will exist because $W$ is well-generated;
see Subsection~\ref{NC(W)-subsection}.  In other words,
$c$ has an eigenvector $v \in V$ fixed by none of the reflections, and whose $c$-eigenvalue
$\zeta_h$ is a primitive $h^{th}$ root of unity.

Let 
\begin{equation}
\label{NC(W)-first-definition}
NC(W):=\{w \in W: \ell_R(w) + \ell_R(w^{-1} c) = n\}
\end{equation}
where $\ell_R$ is the {\it absolute length function}  
\begin{equation}
\label{absolute-length-definition}
\ell_R(w):=\min\{\ell: w = r_1 r_2 \cdots r_\ell \text{ for some }r_i \in R\}.
\end{equation}
The initials ``$NC$'' in $NC(W)$ are motivated by the
special case where $W$ is the Weyl group of type $A_{n-1}$, and
the set $NC(W)$ bijects with the {\it \underline{n}on\underline{c}rossing partitions} 
of $\{1,2,\ldots,n\}$;
see Subsection~\ref{CSP-review} below.

The fact that $R$ is stable under $W$-conjugacy implies that
the cyclic group $C:=\langle c \rangle$ acts on the set $NC(W)$ by conjugation: $w \in NC(W)$ implies
$cwc^{-1} \in NC(W)$.

\begin{theorem}
\label{main-result}
In the above setting, the triple $(X,X(q),C)$ given by
$$
\begin{aligned}
X&=NC(W) \\
X(q)&=\Cat(W,q)\\
C&= \langle c \rangle \cong \ZZ/h\ZZ
\end{aligned}
$$
exhibits the cyclic sieving phenomenon defined in \cite{RSW}:  for any element $c^i$ in $C$,
the subset $X^{c^i}$ of elements in $X$ fixed by $c^i$ has cardinality
$|X^{c^i}|  = [ X(q) ]_{q =\zeta_h^i}$.
\end{theorem}

After explaining a bit more of the background and definitions in Section~\ref{further-background-section},
Theorem~\ref{main-result} is proven in Section~\ref{proof-section}.  
The (partly conjectural) interpretation of $\Cat(W,q)$ is discussed in 
Section~\ref{interpretation-section}, leading to speculation on more conceptual proofs of 
Theorem~\ref{main-result} in Section~\ref{speculation-section}.  Section~\ref{variations-section} remarks
on some conjectural variations.

\tableofcontents

\section{Further background and definitions}
\label{further-background-section}

\subsection{The cyclic sieving phenomenon}
\label{CSP-review}
 
In \cite{RSW}, the following enumerative phenomenon was defined, and then
identified in many examples.  

\begin{definition}
\label{CSP-definition}
Let $C$ be a finite cyclic group, $X$ a finite
$C$-set, and
$X(q)$ a polynomial in $\ZZ[q]$.
Say that the
triple $(X,X(q),C)$ exhibits the {\it cyclic sieving phenomenon} (abbreviated $CSP$
hereafter) if
for every $c \in C$ and any $\zeta \in \CC^\times$ a root of unity having
the same multiplicative order as $c$,
one has
$$
\left[ X(q) \right]_{q=\zeta} = |X^c|
$$
where $X^c:=\{x \in X: c(x)=x\}$.  In particular, this requires $X(1) = |X|$, so
that one thinks of $X(q)$ as a generating function in $q$ for the set $X$;
in many of the observed instances, $X(q)$ has nonnegative coefficients
and is of the form $X(q)=\sum_{x \in X}q^{s(x)}$ for some statistic $s(x)$ defined
for $x \in X$.
\end{definition}

Theorem~\ref{main-result} generalizes the following instance of a CSP from \cite{RSW}.  
Given the $n$-element set $[n]:=\{1,2,\ldots,n\}$,
think of it as labelling the vertices of a convex $n$-gon in a circular order.
Then define a {\it noncrossing partition} of $[n]$ to be a disjoint decomposition of $[n]=\bigsqcup_i B_i$ 
such that different blocks $B_i \neq B_j$ label sets of vertices whose
convex hulls are pairwise disjoint.
Equivalently, for $1 \leq a < b < c < d \leq n$ one never has $a,c$ together in a block $B_i$
and $b,d$ together in a different block $B_j$.
Let $X=NC(n)$ be the set of all noncrossing partitions of $[n]$,
with $C=\ZZ/n\ZZ$ acting by cyclic rotations of $[n]$.  Let $X(q)$ be the following
$q$-{\it Catalan number}
$$
X(q):= C_n(q) := \frac{1}{[n+1]_q} \qbin{2n}{n}
$$
introduced originally by F\"urlinger and Hofbauer \cite{FurlingerHofbauer}.
Here 
$$
\begin{aligned}
\qbin{n}{k}& := \frac{[n]!_q}{[n-k]!_q [k]!_q}, \\
[n]!_q &:= [n]_q [n-1]_q \cdots [2]_q [1]_q
\end{aligned}
$$
with $[n]_q=\frac{q^n-1}{q-1}$ as before.
In \cite{RSW}, the following result was proven, essentially by brute force calculation.

\begin{theorem}\cite[Theorem 7.2]{RSW}
\label{RSW-result}
The triple $(X,X(q),C)$ given by
$$
\begin{aligned}
X&=NC(n) \\
X(q)&=C_n(q)\\
C&= \ZZ/n\ZZ
\end{aligned}
$$
exhibits the $CSP$.
\end{theorem}

To view this as an instance of Theorem~\ref{main-result}, recall that
$W=\symm_n$ the symmetric group on $n$ letters acts on $\CC^n$ by permuting coordinates,
and is generated by the set $R$ of all transpositions, which act as reflections.  
In fact, $W$ is the Weyl group of type $A_{n-1}$ generated by the $n-1$ adjacent transpositions,
and restricts to a well-generated irreducible reflection group acting on the subspace
$V \cong \CC^{n-1}$ where the sum of coordinates is $0$.
When $W$ acts on $V$, its fundamental degrees are $2,3,\ldots,n-1,n$,
since, for example, the elementary symmetric functions $e_2,\ldots,e_n$
generate $S^W$.  Thus the Coxeter number $h=n$ and one can check that $\Cat(W,q)=C_n(q)$.  
One may choose as a Coxeter element the $n$-cycle $c=(1 \,\, 2 \,\, 3 \cdots n-1 \,\, n)$, since
$c(v)=\zeta_n v$ for the vector $v=(1,\zeta_n,\zeta_n^2,\ldots,\zeta_n^{n-1})$ which is fixed
by no reflections (=transpositions).  One then has a map $NC(W) \rightarrow NC(n)$ which
sends a permutation $w$ in $NC(W)$ to the partition $[n]=\bigsqcup_i B_i$ in which the 
$B_i$ are the cycles of $w$.  This maps turns out to be a bijection, sending the
conjugation action of $C$ to cyclic rotations of $[n]$, so that the
the triples $(X,X(q),C)$ in Theorems~\ref{main-result} and \ref{RSW-result} agree.

\subsection{Degrees, codegrees, well-generation}
\label{well-generation-subsection}

We recall here some of the numerology of degrees and codegrees for complex reflection
groups, including an important feature of well-generated groups.

Recall the {\it degrees} $d_1 \leq \cdots \leq d_n$ of a complex reflection group $W$ are
the degrees of homogeneous polynomial generators for the invariant ring $S^W=\CC[f_1,\ldots,f_n]$, 
where $S:=\Sym(V^*)=\CC[x_1,\ldots,x_n]$.  It turns out that they can defined in
an equivalent way using the {\it coinvariant algebra} $S/(\ff)$,
where $(\ff):=(f_1,\ldots,f_n)=S^W_+$ is the ideal of $S$ generated by the
positive degree $W$-invariants.  Both Chevalley \cite{Chevalley} and 
Shephard and Todd \cite{ShephardTodd}) showed that
$S/(\ff)$ carries the regular representation of $W$.  Consequently
$S/(\ff)$ contains exactly $r$ copies of any irreducible $W$-representation $U$ of dimension $r$.  
The {\it U-exponents} $e_1(U),\ldots,e_r(U)$ are the degrees of the homogeneous components
$S/(\ff)$ in which these $r$ copies of $U$ occur.  It is known (see e.g. \cite[\S2]{OrlikSolomon})
that one can uniquely determine/characterize the degrees for $W$ by saying that
the $V$-exponents are $d_1-1,d_2-1,\ldots,d_n-1$.  One also can define the 
{\it codegrees} $d_1^* \geq \cdots \geq d_n^*$ for $W$ uniquely by saying that
the $V^*$-exponents are $d_1^*+1,d_2^*+1,\ldots,d_n^*+1$.

Part of the numerology of complex reflection groups, first observed by
Orlik and Solomon \cite{OrlikSolomon} via the Shephard and Todd 
classification \cite{ShephardTodd}, is that well-generated groups may be
characterized by a ``duality'' between degrees and codegrees:

\begin{theorem}
\label{duality-equals-well-generated}
Let $W$ be an irreducible complex reflection group. The following assertions are
equivalent:
\begin{itemize}
\item[(i)] $W$ is well-generated,
\item[(ii)] $W$ is a {\it duality group} in the sense
that its degrees $d_1 \leq \cdots \leq d_n$
and codegrees $d_1^* \geq \cdots \geq d_n^*$ satisfy
\begin{equation}
\label{well-generated-degs-codegs}
d_i + d_i^* = d_n.
\end{equation}
\end{itemize}
\end{theorem}

\indent
When $W$ is well-generated, $h:=d_n$ is called the {\it Coxeter number}.
Say that an element $c \in W$ is regular if it has an eigenvector lying in the
complement $V^{reg}$ of the reflecting hyperplanes for $W$.  Say that $c$ is
{\it $\zeta$-regular} if this eigenvector may be taken to have eigenvalue $\zeta$;
when this occurs, say that the multiplicative order $d$ of $\zeta$ is a
{\it regular number} for $W$.  Lehrer and Springer first observed the
following fact using the Shephard and Todd classification;  it was later
proven uniformly by Lehrer and Michel \cite{LehrerMichel}.

\begin{theorem}
\label{regular-numbers}
For any complex reflection group, a positive integer $d$ is a regular
number if and only if $d$ divides as many degrees as it divides codegrees.
\end{theorem}

Note that one always has $d_1^*=0$, and hence a regular number $d$ must be
a divisor of at least one of the degrees $d_i$.  Combining Theorem~\ref{regular-numbers}
with \eqref{well-generated-degs-codegs} shows that for $W$ well-generated,
the Coxeter number $h$ is always a regular number, so there exists a regular
element $c$ with eigenvalue $\zeta_h$, a primitive $h^{th}$ root of unity.
A result of Springer \cite[Theorem 4.2]{Springer} shows that for any $\zeta$, 
all $\zeta$-regular elements are $W$-conjugate, and hence $c$ is unique up 
to conjugacy;  elements of $W$ in this conjugacy class are called {\it Coxeter elements}.

\subsection{The $W$-noncrossing partitions}
\label{NC(W)-subsection}

It was recently observed by the first author \cite{Bessis1, Bessis2} that the foregoing
considerations, along with results of Brady and Watt \cite{BradyWatt} allow one to generalize
noncrossing partitions to all irreducible well-generated
complex reflection groups $W$.  Nice background surveys on these topics and
their relation to other aspects of {\it Catalan combinatorics for reflection groups}
are given by Armstrong \cite{Armstrong1} and Fomin and Reading \cite{FominReading-ParkCity}.

One can use the absolute length function
from \eqref{absolute-length-definition} to define a partial order on $W$ as follows:
$w \leq w'$ if $\ell_R(w) + \ell_R(w^{-1}w') = \ell_R(w')$.
This partial order on $W$ has the identity element $e$ as
minimum element, and any Coxeter element $c$ will be among its maximal elements.
The {\it poset of noncrossing partitions for $W$}, denoted
$NC(W)$, is defined to be the interval $[e,c]_W$ in this partial order; this turns out to
agree with the earlier definition \eqref{NC(W)-first-definition}.
Since the set $R$ is closed under conjugation, the $W$-action on $W$ by
conjugation preserves $\ell_R$, and hence respects the partial order $\leq$ on $W$.
Furthermore, as noted earlier, $C=\langle c \rangle$ acts on $NC(W)$ by conjugation.  
Since $c$ is uniquely defined
up to $W$-conjugacy, the isomorphism type of $NC(W)$ and even its $C$-action
is uniquely defined up to isomorphism by the structure of $W$ as a reflection group.

It turns out that for well-generated groups $W$, the poset
$NC(W)$ shares many of the well-known properties enjoyed by its special instance
in type $A_{n-1}$, the poset $NC(n)$:

\begin{theorem}(Bessis \cite{Bessis2})
\label{NC(W)-properties}
For a well-generated complex reflection group $W$, one has that
\begin{enumerate}
\item[(i)]  $NC(W)$ is a bounded, graded poset, which is self-dual and
locally self-dual,
\item[(ii)] $NC(W)$ is a lattice, and
\item[(iii)] its cardinality $|NC(W)|$ is given by the {\it $W$-Catalan number}
$$
\Cat(W):= \prod_{i=1}^n \frac{ h + d_i }{d_i} 
            \quad \left( = \left[ \Cat(W,q)\right]_{q=1} \right).
$$
\end{enumerate}
\end{theorem}

\noindent
We should remark that the properties listed in part (i) of the theorem
can be deduced in a case-free way, but
properties (ii), (iii) have only been proven so far for well-generated
$W$ by appeal to Shephard and Todd's classification of complex 
reflection groups \cite{ShephardTodd}.

\subsection{$W$-$q$-Catalan numbers}

It is not at all clear that the $W$-$q$-Catalan number defined by
\eqref{W-q-Catalan-definition}
should be a polynomial in $q$, nor that it should have
nonnegative coefficients;  the same is true of their
{\it $m$-extended} or {\it Fuss-Catalan} generalizations defined in
\eqref{m-W-q-Cat-definition} below.

For real reflection groups, these $q$-Catalan and $q$-Fuss-Catalan numbers 
for $W$ seem to have first appeared 
in the paper of Berest, Etingof and Ginzburg \cite{BerestEtingofGinzburg}; 
see Section~\ref{interpretation-section} below.  Their
work implies that for such groups they are polynomials in $q$ with nonnegative coefficients, a fact
which had also been conjectured independently around the same time
by Athanasiadis and Garoufallidis; see \cite[Problem 2.1]{Armstrong2}.
For well-generated complex reflection groups, this still holds, but at the moment
has only been verified by appeal to the classification-- see Section~\ref{interpretation-section}
for a conjectural algebraic interpretation.

\subsection{On the $CSP$ for non-faithful $C$-actions}

It can happen that a triple $(X,X(q),C)$ exhibiting the CSP has the group
$C$ acting non-faithfully on $X$.  Assuming $C \cong \ZZ/h\ZZ$ acts with
kernel of order $d$, its quotient group $C'$ of order $\frac{h}{d}$ acts faithfully.
Then one can check that the CSP for $(X,X(q),C)$ is equivalent to the assertion
that $X(q)=X'(q^d)$ for some polynomial $X'(q) \in \ZZ[q]$, and the triple
$(X,X'(q),C')$ exhibits the CSP.

Although this does not happen in Theorem~\ref{main-result} for type $A$ (that is,
Theorem~\ref{RSW-result}), it
can happen in general, as the conjugacy action of $\langle c \rangle$ on $NC(W)$ may not
be faithful: $c$ may have non-trivial central powers.
The following lemma characterizing this situation is
well-known (the first statement is \cite[Corollary 3.3]{Springer}; the second
statement is an easy consequence).

\begin{lemma}
\label{centerlemma}
For $W$ an irreducible complex reflection group with degrees
$d_1,\dots,d_n$, the center of $W$ is a cyclic group of order
$d:=\gcd(d_1,\dots,d_n)$.  

Furthermore, if $W$ is well-generated, with Coxeter number $h$, and $c$ any Coxeter element
as defined above, then the center of $W$ is generated by $c^\frac{h}{d}$.
\end{lemma}

\section{Proof of Theorem \ref{main-result}}
\label{proof-section}

The proof of Theorem~\ref{main-result} has two pieces.  One is further numerology
about degrees, codegrees, and Coxeter elements, given in Lemma~\ref{centralizer-inheritance} below.
The other is about the structure of $NC(W)$, given in Lemma~\ref{conjectured-bijection} below.

\subsection{A lemma on further numerology}

\begin{lemma}
\label{centralizer-inheritance}
Let $W$ be an irreducible complex reflection group acting on $V$ which
\begin{enumerate}
\item[(a)] is well-generated,
\item[(b)] has Coxeter number $h$,
\item[(c)] admits $c$ as one of its Coxeter elements.
\end{enumerate}

\noindent
Given $d$ dividing $h$, define the subgroup 
$$
W':=\Cent_W(c^{\frac{h}{d}}).
$$

Then $W'$ acts on the $\zeta_h^{\frac{h}{d}}$-eigenspace $V'$ for
$c^{\frac{h}{d}}$ as
a complex reflection group, and inherits all three properties 
(a),(b),(c) verbatim.  

Furthermore, if
$W$ has degrees $\{d_i\}_{i=1}^n$ (resp. codegrees $\{d_i^*\}_{i=1}^n$)
then $W'$ acts on $V'$ with degrees (resp. codegrees) equal to the
submultiset of the $\{d_i\}$ (resp. $\{d_i^*\}$) which are divisible by $d$.
\end{lemma}

\begin{proof}
Springer already proved in \cite[Theorem 4.2]{Springer} that
$W'$ is a complex reflection group acting on $V'$, and
that its degrees are as stated in the last assertion of the theorem.  
The statement about codegrees was later proved by Lehrer and Springer,
\cite[Theorem C]{LS}.

This implies that $W'$ is a duality group, thus that it is well-generated
by Theorem~\ref{duality-equals-well-generated}.
Property (b) for $W'$ is obvious.
For (c), note that $c$ does centralize $c^{\frac{h}{d}}$, so
that it lies in $W'$.  Also if $v \in V^{reg}$ is a regular $\zeta_h$-eigenvector for $c$ then 
$v$ is also a $\zeta_h^{\frac{h}{d}}$-eigenvector for $c^{\frac{h}{d}}$, so $v$ lies
in $V'$ and gives a regular $\zeta_h$-eigenvector for $c$ there.
\end{proof}

\subsection{Body of the proof and the structural lemma}

The proof of Theorem~\ref{main-result} will use
a root-of-unity fact arising frequently in cyclic sieving phenomenona, 
and which is an easy consequence of L'H\^opital's rule.

\begin{proposition}
\label{numerical-prop}
Assume that $m', m$ are positive integers with $m' \equiv m \mod d$,
and let $\zeta$ be a $d^{th}$ root of unity.  Then
$$
\lim_{q \rightarrow \zeta} \frac{[m']_q}{[m]_q} =
\begin{cases}
\frac{m'}{m} & \text{ if }m' \equiv m \equiv 0 \mod d \\
1            & \text{ otherwise.}
\end{cases}
$$
\end{proposition}

To prove the theorem, assume one is given an element $c^i \in C$, say of order $d$.
One must show that
\begin{equation}
\label{desired-equation}
\left[ \Cat(W,q) \right]_{q=\zeta_d}
=
|\{ w \in NC(W) : c^i w c^{-i} = w \}| 
\end{equation}
where $\zeta_d = \zeta_h^i$ is a primitive $d^{th}$ root of unity.

Define the centralizer
 subgroup $W':=\Cent_W(c^i)=\Cent_W(c^{\frac{h}{d}})$, where the
second equality follows from the fact that $c^i, c^{\frac{h}{d}}$ have the same order
$d$, and hence generate the same cyclic subgroup of $C$.

Working on the left side of \eqref{desired-equation}, since $d$ divides $h$, 
one has $d_i + h \equiv d_i \mod d$ for all $i$.  Hence 
$$
\left[ \Cat(W,q) \right]_{q=\zeta_d} = \prod_{\substack{i=1,\ldots,n \\ d | d_i}} \frac{h+d_i}{d_i} = |NC(W')|
$$
where the first equality uses Proposition~\ref{numerical-prop},
and the second uses Lemma~\ref{centralizer-inheritance} along with 
Theorem~\ref{NC(W)-properties} (iii) for $W'$.

Meanwhile, the right side of \eqref{desired-equation} is, by definition,
$|NC(W) \cap W'|$.  Thus Theorem~\ref{main-result} will follow from this
key structural lemma, proven in the next two subsections:

\begin{lemma}
\label{conjectured-bijection}
With the above notation and hypotheses, one has 
$$NC(W) \cap W' = NC(W').$$
\end{lemma}

\noindent
Note one slightly tricky point:  we are regarding $W'$ as a subgroup of
$W$ when we write $NC(W) \cap W'$,
but the definition of $NC(W')$ makes sense when $W'$ is viewed as
a complex reflection group acting on $V'$.
In particular, the subset $R'$ of $W'$ consisting of elements 
acting as reflections on $V'$ is not in general included in $R$, and
the reflection length function $\ell_{R'}$ defined on $W'$ is not
the restriction of $\ell_R$.

The forward inclusion $NC(W)\cap W' \subseteq NC(W')$ in Lemma~\ref{conjectured-bijection}
is straightforward and proven in the next subsection,
using the elementary
linear algebraic interpretation of the partial order on $W$ due to 
Brady and Watt \cite{BradyWatt}.
The reverse inclusion appears harder,
and is proven in the following subsection.  It uses the geometric reinterpretation of $NC(W)$ that
plays an important role in work of the first author \cite{Bessis2}, proving that
complements of complex reflection arrangements are $K(\pi,1)$'s.  It would be nice to have a simpler
proof of this reverse inclusion.

\subsection{The straightforward inclusion}
\label{straightforward-inclusion-subsection}

To prove the inclusion 
\begin{equation}
\label{straightforward-inclusion}
NC(W)\cap W' \subseteq NC(W')
\end{equation}
we recall Brady and Watt's characterization  \cite{BradyWatt} of
the partial order on $W$ via fixed spaces:
$x \leq z$ if and only if the factorization $xy=z$ (that is, where $y:=x^{-1}z$) satisfies
\begin{equation}
\label{Brady-Watt-conditions}
\begin{aligned}
V^x \cap V^y &= V^z,\\
V^x + V^y &=V.
\end{aligned}
\end{equation}

Now given $x \in NC(W) \cap W'$, we can factor $xy=c$ in $W$ 
satisfying conditions \eqref{Brady-Watt-conditions} with $z=c$.  Since $W$ was assumed to
act irreducibly on $V$, we
claim that $V^c=0$.  To see this, note first that the eigenvalues of $c$ are $\{\zeta_h^{1-d_i}\}_{i=1}^n$,
according to Springer \cite[Theorem 4.2]{Springer},
so if $V^c \neq 0$, then there exists some $i$ for which $\zeta_h^{1-d_i}=1$.  As
the degrees $d_i$ of the fundamental $W$-invariants all lie the range $1 \leq d_i \leq h$,
this forces $d_i=1$ for some $i$, and hence the existence of a $W$-invariant linear
form on $V$.  But then the kernel of this $W$-invariant linear form would be
a proper $W$-invariant linear subspace, contradicting irreducibility.

Thus there is a direct sum decomposition
\begin{equation}
\label{fixed-space-direct-sum}
V^x \oplus V^y = V.
\end{equation}

\noindent
Note that $c \in W'$ since $c$ centralizes $c^{\frac{h}{d}}$, and
hence since $x \in W'$, we conclude that $y (=x^{-1}c)$ also lies in $W'$.  Therefore 
$$
\left( x|_{V'} \right) \left( y|_{V'} \right) =c|_{V'}
$$
is a factorization in $W'$.  
Furthermore, intersecting both sides of 
\eqref{fixed-space-direct-sum} with $V'$ gives
\begin{equation}
\label{intersected-direct-sum}
(V')^x \oplus (V')^y = V'.
\end{equation}
Since $(V')^c \subset V^c=0$, the direct sum
\eqref{intersected-direct-sum} shows
$$
\begin{aligned}
(V')^x + (V')^y &= V'\\
(V')^x \cap (V')^y & = 0=(V')^c.
\end{aligned}
$$
Hence $x|_{V'} \in NC(W')$ as desired.

\subsection{The less straightforward inclusion}
\label{harder-inclusion-subsection}

The reverse inclusion to \eqref{straightforward-inclusion} 
is better understood using the geometry of
the quotient space of $V$ under the action of $W$, 
where elements of $NC(W)$ and $NC(W')$ have convenient interpretations.

Let us start with some standard facts from invariant theory (see
for example \cite{Springer}). 
After choosing a system of homogeneous basic invariants $f_1,\dots,f_n$ for $W$
for which $S^W=\CC[f_1,\ldots,f_n]$, one can view $(f_1,\dots,f_n)$
as coordinate functions on the quotient space $V/W$: they realize
an explicit isomorphism $V/W\simeq \mathbf{C}^n$.  

Given a positive integer $d$, let 
$I:= \{ i \in \{1,\dots,n\}: d \text{ divides }d_i\}$
and denote its complement $J:=\{1,\dots,n\} \setminus I$.
A standard consequence of \cite[3.2 and 4.2]{Springer} is that the
quotient space $V'/W'$ gets identified with the linear subspace
of $V/W$ defined by the equations $\{ f_j=0 \}_{j \in J}.$
The {\it discriminant hypersurface} $\mathcal{H}$ (resp. $\mathcal{H}'$)
is the algebraic hypersurface which is the image under the quotient mapping $V \rightarrow V/W$ (resp. $V'/W'$)
of the union of all reflecting hyperplanes for $W$ inside $V$ (resp. for $W'$
inside $V'$).  The braid group $B$ (resp. $B'$) associated with $W$
(resp. $W'$) is defined as the fundamental
group of the complement of $\mathcal{H}$ in $V/W$ (resp. the complement of $\mathcal{H}'$ in $V'/W'$).
It was proven independently by Lehrer and Denef-Loeser 
(see \cite[Theorem 1.9 (3)]{Bessis2}) that 

\begin{equation}
\label{Denef-Loeser-Lehrer}
\mathcal{H}\cap (V'/W')= \mathcal{H}'.
\end{equation}
This gives an inclusion $V'/W'-\mathcal{H}' \hookrightarrow
V/W -\mathcal{H}$, and hence a natural morphism $B'\rightarrow B$. 
Combined with the natural epimorphisms from the braid groups to the
associated reflection groups, we get a commutative diagram

\begin{equation}
\label{braid-group-diagram}
\xymatrix{ B' \ar[r]  \ar@{>>}[d]_{\pi'}  & B  \ar@{>>}[d]^{\pi} \\
 W'  \ar@{^{(}->}[r] & W . }
\end{equation}

We now explain how $NC(W)$ and $NC(W')$ are 
related to all this, assuming that $W$ is well-generated and $d$ divides $h:=d_n$.
In \cite[section 6]{Bessis2}, one considers certain paths in $V/W-\mathcal{H}$
called \emph{tunnels}. For a given
path $\gamma:[0,1] \rightarrow V/W -\mathcal{H}$,
being a tunnel essentially amounts to saying that, along the path $\gamma$,
\begin{enumerate}
\item[$\bullet$] the coordinates $(f_1,\dots,f_{n-1})$ are
constant, and
\item[$\bullet$] the imaginary part $\im f_n$ is constant, but
\item[$\bullet$] the real part $\re f_n$ is increasing.
\end{enumerate}
The similar notion for $W'$ involves the coordinates $(f_i)_{i\in I}$
instead of $(f_1,\dots,f_n)$ -- the coordinate whose real part is allowed to
increase along a tunnel for $W'$ still is the $f_n$-coordinate.
(Note that since $W$ is well-generated, and $d$ is a divisor of
$h:=d_n$, one indeed knows that $n$ lies in $I$).
From this description, one immediately concludes that
if a path $\gamma':[0,1] \rightarrow V'/W' -\mathcal{H}'$ is a tunnel
(with respect to $W'$), then, as a path in $V/W - \mathcal{H}$ (which
contains $V'/W' -\mathcal{H}'$), $\gamma'$ is also a tunnel (with respect
to $W$).

Another feature of tunnels is that, to any tunnel $\gamma$ with respect
to $W$, one may associate an element $b_{\gamma}$ in $B$, the fundamental group
of $V/W -\mathcal{H}$.  The definition a tunnel implies that
its endpoints lie inside a certain dense open subset ${\mathcal U}$
of $V/W -\mathcal{H}$, defined in \cite[Defintion 6.1]{Bessis2} by the following condition:
a point $p$ in $V/W$ lies in ${\mathcal U}$ if and only if there exist no 
point $q$ on the hypersurface $\mathcal H$ for which
\begin{equation}
\label{basepoint-definition}
\begin{aligned} 
f_i(p)&=f_i(q) \text{ for }i = 1,2,\ldots, n-1, \text{ and } \\
\re f_n(p)&=\re f_n(q), \text{ and } \\
\im f_n(p)& \leq \im f_n(q). \\
\end{aligned}
\end{equation}
It turns out \cite[Lemma 6.2]{Bessis2} that ${\mathcal U}$ is contractible, and 
therefore can be used as the ``basepoint'' for defining the braid group $B$ 
in \eqref{braid-group-diagram} as the (relative) fundamental group 
$$
B:=\pi_1(V/W-\mathcal{H}, \,\, {\mathcal U}).
$$ 
One uses the similarly-defined contractible base point $\mathcal U'$ inside $V'/W'$
when defining the braid group $B'$ in \eqref{braid-group-diagram}.
The horizontal map $B \rightarrow B'$ is well-defined because these 
contractible base points are compatible:  the definition \eqref{basepoint-definition}
together with \eqref{Denef-Loeser-Lehrer} shows that 
$$
{\mathcal U}'=\left(V'/W'-\mathcal{H}'\right) \cap {\mathcal U}.
$$

Elements $b_{\gamma}$ of $B$ obtained from tunnels $\gamma$
are said to be \emph{simple}.  Similarly, tunnels with respect to $W'$ give rise to \emph{simple}
elements in $B'$.  
Temporarily denote by $S$ (resp. $S'$) the set of simple elements
in $B$ (resp. $B'$).
Proposition 8.9 in \cite{Bessis2} asserts that $\pi$ (resp. $\pi'$)
restricts to a bijection between $S$ and $NC(W)$ (resp. between
$S'$ and $NC(W')$).

The desired inclusion $NC(W') \subseteq NC(W)\cap W'$ is now straightforward
to deduce. Let $w'\in NC(W')$, so that $w' \in W'$ trivially. There
exists a simple element $s'\in \pi^{-1}(w')$, represented
by a tunnel $\gamma'$ in $V'/W' -\mathcal{H}'$. But since $\gamma'$
is also a tunnel in $V/W -\mathcal{H}$, it defines a simple
element $s\in B$. All these constructions are compatible with the
the diagram \eqref{braid-group-diagram} and prove that $w'=\pi'(s') = \pi(s)$
lies in $NC(W)$.

This concludes the proof of the reverse inclusion in Lemma~\ref{conjectured-bijection}, and hence
the proof of Theorem~\ref{main-result}.

\section{Interpreting $\Cat(W,q)$}
\label{interpretation-section}

  Here we review the interpretation of $\Cat(W,q)$ for real reflection groups
due to Berest, Etingof and Ginzburg \cite[Theorem 1.10(iii)]{BerestEtingofGinzburg}.  We also
conjecture how this might generalize to well-generated complex reflection groups $W$, as well as
to some relatives and generalizations of $\Cat(W,q)$.  A crucial idea appears originally in the
work of Haiman \cite[\S2.5]{Haiman} on diagonal harmonics, where the interpretation of
many formulae relied on the existence of a system of parameters carrying the following extra structure.

\begin{definition}
Let $W$ be a complex reflection group acting on $V=\CC^n$, and $S:=\Sym(V^*)$ the
algebra of polynomial functions on $V$.  Given an $n$-dimensional representation $U$ of $W$,
say that a collection $\thetas=\{\theta_1,\ldots,\theta_n\}$
of $n$ elements in $S$ which are all homogeneous of the same degree $p$ form
a {\it homogeneous system of parameters (hsop) carrying $U(-p)$} if 
\begin{enumerate}
\item[$\bullet$] they are a system of parameters for $S$, meaning that they are
algebraically independent and the quotient $S/(\thetas)$ is finite dimensional over $\CC$, and
\item[$\bullet$] their span $\CC\theta_1 + \cdots + \CC\theta_n$ is a $W$-stable subspace
inside the $p^{th}$ homogeneous component $S_p$, and carries a $W$-representation equivalent to $U$.
\end{enumerate}
\end{definition}

Finding such hsop's can be difficult.  For real reflection groups $W$,
the work of Berest, Etingof and Ginzburg finds such hsop's 
using subtle results (including work of Gordon \cite{Gordon}) on the representation theory
of {\it rational Cherednik algebras}-- the authors are very grateful to S. Griffeth,
for pointing them to these results.  Once one has such hsop's, various character and Hilbert series
formulae will ensue by fairly standard calculations.  

The following proposition collects several of these formulae,
some of which may be found stated for the case of real reflection groups in 
Berest, Etingof and Ginzburg \cite[Theorems 1.6, 1.7, 1.10]{BerestEtingofGinzburg}
and Gordon \cite[\S 5]{Gordon}.  For notational purposes, recall that if $\sigma$ lies
in the Galois group $\Gal(\QQ[\zeta_{|W|}]/\QQ)$, one can form 
the {\it Galois conjugate representation} 
$U^\sigma$ in which an element $w \in W$ acting on $U$ by a given matrix $A(w)$ will
act by the matrix $A(w)^\sigma$ on $U^\sigma$.  In particular, since $W$ may be taken to act on 
$U$ by unitary transformations, when $\sigma$ is
complex conjugation, the Galois conjugate $U^\sigma$ is equivalent to the representation $U^*$ 
{\it contragredient} (or {\it dual}) to $U$.

\begin{proposition}
\label{many-formula-proposition}
Let $W$ be a complex reflection group acting on an $n$-dimensional complex vector
space $V$, and $U$ an $n$-dimensional representation of $W$.  If $\thetas$ in $S$ is an
hsop carrying $U(-p)$, then the $W$-representation $S/(\thetas)$ has Hilbert series
\begin{equation}
\label{Hilb-of-S-mod-theta}
\Hilb(S/(\thetas),q) = \left( [p]_q \right)^n.
\end{equation}

If furthermore $U=V^\sigma$ is a Galois conjugate of the representation $V$, then the character of
$W$ acting on $S/(\thetas)$ is
\begin{equation}
\label{general-character-formula}
\chi_{S/(\thetas)}(w)= p^{\dim_\CC(V^w)} \prod_{w-\text{eigenvalues }\lambda \neq 1} 
                            \frac{1-\sigma(\lambda)}{1-\lambda},
\end{equation}
and the $W$-invariant subspace of $S/(\thetas)$ has Hilbert series
\begin{equation}
\label{general-product-formula}
\Hilb(\left( S/(\thetas) \right)^W , q) = \prod_{i=1}^n \frac{[p+e_i(U^*)]_q}{[d_i]_q}
\end{equation}
where $e_1(U^*),\ldots,e_n(U^*)$ denote the $U^*$-exponents defined
in Subsection~\ref{well-generation-subsection}.

In particular, when $U=V^*$ is the contragredient representation to $V$, one has
\begin{equation}
\label{V*-formulae}
\begin{aligned}
\chi_{S/(\thetas)}(w)&= p^{\dim_\CC(V^w)} \\
\Hilb((S/(\thetas))^W ,q) &= \prod_{i=1}^n \frac{[p+d_i-1]_q}{[d_i]_q},
\end{aligned}
\end{equation}
and when $U=V$, one has
\begin{equation}
\label{V-formulae}
\begin{aligned}
\chi_{S/(\thetas)}(w)&= p^{\dim_\CC(V^w)} (-1)^{n-\dim_\CC V^w} \overline{\det(w)} \\
\Hilb((S/(\thetas))^W,q) &= \prod_{i=1}^n \frac{[p+d_i^*+1]_q}{[d_i]_q}
\end{aligned}
\end{equation}
\end{proposition}

\begin{proof}
Because $S$ is a Cohen-Macaulay ring, the hsop $\thetas$ is also an $S$-regular
sequence. Consequently, the Koszul complex
for $\thetas$ provides a finite free resolution for $S/(\thetas)$ as an $S$-module:
\begin{equation}
\label{Koszul-resolution}
0 \rightarrow S \otimes \wedge^n U(-p)  \rightarrow 
   \cdots \rightarrow
   S \otimes \wedge^j U(-p)  \rightarrow 
    \cdots \rightarrow S \rightarrow S/(\thetas) \rightarrow 0.
\end{equation}
Here all tensor products are over $\CC$, and $U(-p)$ denotes a 
$\CC$-vector space and representation of $W$ isomorphic to $U$, with its
basis decreed to be of polynomial degree $p$ in order to make the Koszul boundary
maps homogeneous.

Taking Euler characteristics in each homogeneous polynomial degree of
\eqref{Koszul-resolution}, one obtains the following virtual
isomorphism of graded $W$-representations:
\begin{equation}
\label{virtual-isomorphism}
S/(\thetas) \cong \sum_{j=0}^n (-1)^j \,\, S \otimes \wedge^j U(-p).
\end{equation}
Then \eqref{Hilb-of-S-mod-theta} follows by taking Hilbert series of the two sides
of \eqref{virtual-isomorphism}:
$$
\begin{aligned}
\Hilb(S/(\thetas),q) &= \sum_{j=0}^n (-1)^j \Hilb(S \otimes \wedge^j U(-p),q)\\
                     &= \Hilb(S,q) \left[ \sum_{j=0}^n t^j \Hilb(\wedge^j U(-p),q) \right]_{t=-1}\\
                     &= \frac{1}{(1-q)^n} \left[ \sum_{j=0}^n t^j \binom{n}{j} (q^p)^j \right]_{t=-1}\\
                     &= \left[ \left( \frac{1+tq^p}{1-q} \right)^n \right]_{t=-1}\\
                     &= \left( [p]_q \right)^n.
\end{aligned}
$$
Now assume that $U=V^\sigma$, and define the bigraded trace of $w$ 
$$
\begin{aligned}
\chi_{S \otimes \wedge U(-p)}(w; t, q)
   & := \sum_{j=1}^n \sum_{i \geq 0} \chi_{(S \otimes \wedge^j U(-p))_i}(w) \,\, t^j q^i \\
   & = \prod_{w-\text{eigenvalues }\lambda} \frac{1+\sigma(\lambda)t q^p}{1-\lambda q}.
\end{aligned}
$$
Then \eqref{virtual-isomorphism} implies that
the singly graded trace of $w$ on $S/(\thetas)$ will be
$$
\begin{aligned}
\chi_{S/(\thetas)}(w,q) &:=\sum_{ j \geq 0 } \chi_{\left(S/(\thetas)\right)_i}(w) q^i \\
                      & =\left[ \chi_{S \otimes \wedge U(-p)}(t, q) \right]_{t=-1} \\
                      & = \prod_{w-\text{eigenvalues }\lambda} \frac{1-\sigma(\lambda) q^p}{1-\lambda q}.
\end{aligned}
$$
and taking the limit as $q$ approaches $1$ yields \eqref{general-character-formula}.

To prove \eqref{general-product-formula}, note that 
$$
\begin{aligned}
\Hilb((S/(\thetas))^W , q ) &= \sum_{j=0}^n (-1)^j 
                          \Hilb \left( ( S \otimes \wedge^j U(-p)  )^W,q \right) \\
                    &= \left[ \Hilb \left( ( S \otimes \wedge U(-p) )^W , t, q \right) \right]_{t=-1}
\end{aligned}
$$
When $U=V^\sigma$ is a Galois conjugate of $V$, 
Orlik and Solomon~\cite[Theorem 3.1]{OrlikSolomon} analyzed the structure of
$( S \otimes \wedge U )^W$ considered as an $S^W$-module,
generalizing work of Solomon~\cite{Solomon} when $U=V^*$.  Their result says
that $( S \otimes \wedge U )^W$ 
is free over $S^W$, and isomorphic to an exterior algebra over $S^W$.  The
exterior generators $\{h_i \otimes u_i\}_{i=1}^n$ lie in $(S \otimes U)^W$,
and have $h_1,\ldots,h_n$ of polynomial degrees matching the $U^*$-exponents $e_1(U^*),\ldots,e_n(U^*)$.
Therefore 
$$
\begin{aligned}
\Hilb \left( ( S \otimes \wedge U(-p) )^W , t, q \right)
  &= \Hilb(S^W,q) \prod_{i=1}^n (1+tq^{p+e_i(U^*)}) \\
  & =  \prod_{i=1}^n \frac{1+tq^{p+e_i(U^*)}}{1-q^{d_i}}.
\end{aligned}
$$
which yields \eqref{general-product-formula} upon setting $t=-1$.  

The remaining formulae in the theorem now
follow from the characterizations of the degrees $d_i$ and codegrees $d_i^*$ 
in terms of the $V, V^*$-exponents, and from the fact that 
$\frac{1-\bar\lambda}{1-\lambda}=-\overline{\lambda}$
for any root-of-unity $\lambda \neq 1$.
\end{proof}

The interpretation for $\Cat(W,q)$ comes from applying
formula \eqref{V*-formulae} to part (i) of the following conjecture.

\begin{conjecture}
\label{hsop-conjecture}
Let $W$ be an irreducible well-generated finite reflection group, with Coxeter number $h=d_n$.
\begin{enumerate}
\item[(i)] If $p \equiv 1 \, \mod \, h$ then $S$ contains an hsop $\thetas$ carrying $V^*(-p)$.
\item[(ii)] If $p \equiv -1 \, \mod \, h$ then $S$ contains an hsop $\thetas$ carrying $V(-p)$.
\end{enumerate}
\end{conjecture}

\begin{remark}
\label{known-hsops}
We explain here the parts of Conjecture~\ref{hsop-conjecture} known for various classes of
well-generated complex reflection groups.  

Recall that the Shephard and Todd classification of complex reflection groups 
\cite{ShephardTodd} contains only one infinite family:  the groups
$G(d,e,n)$ for $d,e,n$ positive integers with $e$ dividing $d$.
This group $G(d,e,n)$ is the set of $n \times n$ matrices having exactly one
nonzero entry in each row and column, requiring further that this entry must be
a complex $d^{th}$ root of unity, and that the product of these $n$ 
nonzero entries is a $\frac{d}{e}^{th}$ root of unity.
Not all of the groups $G(d,e,n)$ are well-generated-- the well-generated
subfamilies are 
$$
\begin{aligned}
G(1,1,n) & =\text{ the Weyl group of type }A_{n-1},\\
G(d,1,n) & \text{ for }d \geq 2,\\
G(e,e,n) & \text{ for }e \geq 2.
\end{aligned}
$$

Both parts (i) and (ii) of Conjecture~\ref{hsop-conjecture} are easily verified for all
the groups $G(d,e,n)$ with $d \geq 2$, whether or not they are well-generated.
Here $W$ has degrees $d,2d,3d,\ldots,(n-1)d,\frac{nd}{e}$,
so that 
 \begin{equation}
 \label{explicit-classical-h}
 h=\max\left\{(n-1)d,\frac{nd}{e}\right\} =
 \begin{cases}
 (n-1)d       & \text{ if } e>1, \\
  nd           & \text{ if } e=1. 
 \end{cases}
 \end{equation}
Note that $\thetas:=\{x_1^p,\ldots,x_n^p\}$ gives an hsop carrying\footnote{A related
hsop also arises naturally from the viewpoint of rational Cherednik algebras;
see Griffeth \cite[\S 8]{Griffeth}}  $U(-p)$ with
\begin{equation}
\label{G(d,e,n)-hsops}
\begin{cases}
U= V^* &\text{ if }p\equiv +1 \, \mod \, d,\\
U= V   &\text{ if }p \equiv -1 \, \mod \, d.
\end{cases}
\end{equation}
Since $d$ divides $h$ by \eqref{explicit-classical-h},
the assumption that $p \equiv \pm 1\, \mod \, h$ in Conjecture~\ref{hsop-conjecture}
is stronger than assumption $p \equiv \pm 1\, \mod \, d$ in \eqref{G(d,e,n)-hsops}.

Both parts of Conjecture~\ref{hsop-conjecture} hold for the remaining subfamily
$G(1,1,n)$, that is, the Weyl group of type $A_{n-1}$, but the appropriate 
hsop's are not nearly as trivial to construct. Here $W$ only acts irreducibly
after one restricts to the subspace $V \subset \CC^n$ where the sum of coordinates is
zero, and the degrees are $2,3,\ldots,n$ with $h=n$, as explained in Section~\ref{CSP-review}.
Haiman~\cite[Prop. 2.5.3]{Haiman} proves that there is an hsop carrying $V(-p) (\cong V^*(-p))$
if and only if $\gcd(p,n)=1$, via an elementary construction that partly uses an idea of H. Kraft.
In particular, such hsop's exist for $p \equiv \pm 1 \mod h$, as needed in
 Conjecture~\ref{hsop-conjecture}.  

This means that Conjecture~\ref{hsop-conjecture} holds for all of the well-generated
complex reflections within the infinite family $G(d,e,n)$, so there are only finitely
many exceptional cases where it is not known.

In addition, all {\it real} reflection groups are well-generated.  For such groups one has
$V=V^*$, and at least Part (i) of Conjecture~\ref{hsop-conjecture} 
follows from the previously mentioned work of Berest, Etingof and Ginzburg  and of Gordon on
the representation theory of the {\it rational Cherednik algebra} ${\HH}_c$.  This noncommutative
algebra involves the choice of complex parameters $c_s$ for each conjugacy class of
reflections $s \in W$;  for our purposes one must choose all the parameters $c_s$ to
be the same constant $m+\frac{1}{h}$, where $m$ is the unique nonnegative integer satisfying
$p=mh+1$.  A version of the Poincar\'e-Birkhoff-Witt theorem for
this noncommutative algebra $\HH_c$ shows that it contains within it a subalgebra
isomorphic to $S=\Sym(V^*)$ as well as a subalgebra isomorphic to $\CC[W]$.  
For this value of the parameters $c_s$, 
there is an important finite dimensional irreducible ${\HH}_c$-module $L(\triv)$,
whose $S$-module structure turns out to be isomorphic to $S/(\thetas)$.  
The key to seeing this is the crucial and highly nontrivial 
{\it BGG resolution} \cite[Theorems 2.3, 2.4]{BerestEtingofGinzburg}
of $L(\triv)$ by the Verma-like $\HH_c$-modules $M(\wedge^j V^*)$:
\begin{equation}
\label{BGG-resolution}
0 \rightarrow M(\wedge^n V^*) \rightarrow \cdots 
  \rightarrow M(\wedge^2 V^*) \rightarrow M(V^*) 
  \rightarrow M(\triv) \rightarrow L(\triv) \rightarrow 0
\end{equation}
We sketch here why restriction of this resolution to its $S$-module structure 
gives the desired Koszul resolution \eqref{Koszul-resolution} for the case $U=V^*=V$.
Upon restriction to its $S$-module structure, each $M(\wedge^j V^*)$ is a free $S$-module of
rank $\binom{n}{j}$.  The image of $M(V^*)\cong S^n$ will be generated by
$n$ elements $\thetas:=\{\theta_1,\ldots,\theta_n\}$ inside $M(\triv) \cong S$.
One knows, from the $\CC[W]$-module structure of the resolution,
that the span of the $\thetas$ carries the representation $V^*$.  By computing the eigenvalue of
the diagonalizable element ${\mathbf h}$ from \cite{BerestEtingofGinzburg} (whose eigenspaces
characterize the homogeneous components for the grading inside the 
$\HH_c$-modules $M(\wedge^j V^*)$), one concludes that the
$\thetas$ lie in degree $p$, so they carry $V^*(-p)$.  The finite dimensionality of
$L(\triv) \cong S/(\thetas)$ implies $\thetas$ is an hsop.

\begin{question}
Can the representation theory of rational Cherednik algebras be used to prove
Conjecture~\ref{hsop-conjecture} part (ii) also?  

Can it be used to prove both parts (i) and (ii) for the 
well-generated complex reflection groups outside of the real reflection groups?
\end{question}
\end{remark}
\noindent
We mention here that recent work of Berenstein and Burman \cite{BerensteinBurman},
using rational Cherednik algebras to define what they call {\it quasiharmonic polynomials} inside
the polynomial algebra $S$, may be of help in identifying the conjectured
hsop's.

The next proposition summarizes how the two parts of Conjecture~\ref{hsop-conjecture}
interpret the $W$-$q$-Catalan number $\Cat(W,q)$, as well as its
{\it Fuss-Catalan} generalization, and its close relative sometimes called the
{\it positive} $W$-$q$-{\it Catalan number}.
The proof is an immediate consequence of Proposition~\ref{many-formula-proposition}.

\begin{proposition}
\label{Cat-consequences}
When Conjecture~\ref{hsop-conjecture} (i) holds for $W$,
one has the following interpretation for the $q$-analogue of the
{\it $m$-extended} or {\it Fuss-Catalan} numbers:
\begin{equation}
\label{m-W-q-Cat-definition}
\begin{aligned}
\Cat^m(W,q) &:= \prod_{i=1}^{n} \frac{ [m h+d_i]_q }{ [d_i]_q }
= \prod_{i=1}^{n} \frac{ [m h+e_i+1]_q }{ [e_i+1]_q}\\
            &= \Hilb((S/(\thetas))^W,q),
\end{aligned}
\end{equation}
where $\thetas$ is the hsop carrying $V^*(-p)$ for
$p=1+mh$, and the $e_1,\ldots,e_n$ are the $V$-exponents.

When Conjecture~\ref{hsop-conjecture} (ii) holds for $W$,
one has the following interpretation for the $q$-analogue of the ``positive'' version of the
$m$-extended/Fuss-Catalan numbers:
$$
\begin{aligned}
\Cat_+^m(W,q) &:= \prod_{i=1}^{n} \frac{ [m h+d_i^*]_q }{ [d_i]_q }
= \prod_{i=1}^{n} \frac{ [m h+e^*_i-1]_q }{ [e_i+1]_q}\\
            &= \Hilb((S/(\thetas))^W,q)
\end{aligned}
$$
where $\thetas$ is the hsop carrying $V(-p)$ for
$p=-1+mh$, and the $e_1^*,\ldots,e_n^*$ are the $V^*$-exponents.
\end{proposition}

\begin{remark} \rm \ 
\label{positive-Catalan-remark}
When $W$ is a real reflection group (so $V^*=V$ and $\{e^*_i\}_{i=1}^n=\{e_i\}_{i=1}^n$), 
the positive version $\Cat^m_+(W)$ of the $W$-Fuss-Catalan numbers, say for $m=1$, arises in
several places: 
\begin{enumerate}
\item[$\bullet$] in the study of $NC(W)$, as the rank of the {\it homology}
of the proper part of the poset, 
\item[$\bullet$] in the {\it cluster algebras} of finite type, when counting {\it positive
clusters} (that is, those containing only positive roots), and
\item[$\bullet$] in the {\it antichains of positive roots} (for Weyl groups), when counting
those antichains which contain {\it no simple roots}.
\end{enumerate}
We refer the reader to Armstrong \cite[Chapter 1]{Armstrong1} and
Fomin and Reading \cite{FominReading-ParkCity} for more background and history on these topics.
\end{remark}

\section{Speculation}
\label{speculation-section}

  Here we review a method from \cite{RSW},
based on Springer's theory of regular elements, that provides a more
conceptual proof for certain cyclic sieving phenomena.
We then speculate on how this might extend to prove Theorem~\ref{main-result} 
more conceptually.  

We also review the connection, originally due to Haiman \cite[\S7]{Haiman},
between $\Cat(W)$ and the finite torus when $W$ is a Weyl group, and speculate on how this might 
lead to a related CSP.

\subsection{Review of Springer's theory and a known CSP}

Recall that for any complex reflection group $W$,
both Chevalley \cite{Chevalley} and Shephard and Todd \cite{ShephardTodd}
showed the coinvariant algebra $S/(\ff)$ carries the regular representation of $W$.
That is, one has an isomorphism of $\CC[W]$-modules 
\begin{equation}
\label{Chevalley-isomorphism}
S/(\ff) \underset{\CC[W]-\text{mod}}{\cong} \CC[W].
\end{equation}
One of the main results of Springer \cite{Springer} can be rephrased as follows.
Let $C=\langle c \rangle$ be the cyclic group generated by any regular element
$c$ in $W$, with eigenvalue $\zeta$ on some eigenvector $v$ lying in 
$V^{reg}$;  see Subsection~\ref{well-generation-subsection} for the terminology.
Then Springer showed \eqref{Chevalley-isomorphism}
extends to an isomorphism of $\CC[W \times C]$-modules
\begin{equation}
\label{Springer-isomorphism}
S/(\ff) \underset{\CC[W \times C]-\text{mod}}{\cong} \CC[W].
\end{equation}
Here $W, C$ respectively act by {\it left, right} multiplication on $\CC[W]$,
and $W$ acts on $S/(\ff)$ by its usual linear substitutions, while $C$ acts on
$S/(\ff)$ by the scalar substitutions $c(x_i)=\zeta^{-1} x_i$ for all $i$.  This
means that $c$ acts on an element $f$ which is homogeneous of degree $d$ by 
$c(f)=\zeta^{-d}(f)$.

As noted in \cite[\S 8]{RSW}, given an {\it arbitrary} subgroup $W'$ of $W$,
one can one restrict the isomorphism \eqref{Springer-isomorphism}
to the $W'$-fixed subspaces on both sides, yielding a $\CC[C]$-module isomorphism
\begin{equation}
\label{CSP-isomorphism}
\left( S/(\ff) \right)^{W'} \underset{\CC[C]-\text{mod}}{\cong} \CC[W]^{W'} \cong \CC[W'\backslash W]
\end{equation}
in which $C$ acts on $\left( S/(\ff) \right)^{W'}$ by the same scalar substitutions 
related to the grading as before, and $C$ acts on the right cosets $W'w$ via right
multiplication.  Taking the character/trace of any element of $C$ then immediately
yields the following

\begin{theorem}\cite[Theorem 8.2]{RSW}
\label{Springer-CSP}
Let $W$ be a complex reflection group, $W'$ any subgroup of $W$,
and $c$ any regular element of $W$.  Then the triple $(X,X(q),C)$ with
$$
\begin{aligned}
X &=W/W' \\
X(q) &= \Hilb( (S/(\ff))^{W'}, q )
\,\, \left( = \Hilb( S^{W'}/(\ff), q )
            = \frac{\Hilb(S^{W'},q)}{\Hilb(S^W,q)}
\right)\\
C &=\langle c \rangle
\end{aligned}
$$
exhibits the CSP.
\end{theorem}

\subsection{Speculation on NC(W)}

When $W$ is a well-generated complex reflection group,
there is hope of extending Theorem~\ref{Springer-CSP}, beginning
with the following observation.

\begin{proposition}
\label{thetas-inside-fs}
Let $W$ be a well-generated complex reflection group, and assume one has
an hsop $\thetas$ carrying $V^*(-p)$ for $p=h+1$ as in Conjecture~\ref{hsop-conjecture} Part (i).

Then $(\thetas) \subset (\ff)$, so that $S/(\thetas)$ surjects onto 
the coinvariant algebra $S/(\ff)$.
\end{proposition}
\begin{proof}
By definition of the codegrees $d_i^*$, one finds copies of $V^*$ in
$S/(\ff)$ only in degrees $d_n^*+1 \leq \cdots \leq d_1^*+1$
of $S/(\ff)$.  Since $d_i^*=h-d_i$ for well-generated $W$, one concludes
that the highest degree of $S/(\ff)$ containing a copy of $V^*$ is $h-d_1+1 < h+1$.
Thus there is no copy of $V^*$ in degree $p=h+1$ of $S/(\ff)$.  By complete reducibility
of $\CC[W]$-representations, one can decompose $S = H \oplus (\ff)$ where
$H$ is a graded vector space in which $H_p$ is
a $W$-invariant complement to $(\ff)_p$ inside $S_p$ for each $p$.  As $H_p$ is 
$W$-isomorphic to $(S/(\ff))_p$, it contains no copy of $V^*$.  Hence 
$\thetas$ lies in $(\ff)_p$, and $(\thetas) \subset (\ff)$.
\end{proof}

This leads to our first speculation.

\begin{question}
\label{enlarge-NC(W)-question}
In the setting of Proposition~\ref{thetas-inside-fs}, 
does there exist a finite set $P$ 
\begin{enumerate}
\item[$\bullet$] with cardinality $|P|=(h+1)^n$,
\item[$\bullet$] carrying a $W \times C$-action, whose $W$-orbits $W\backslash P$ are 
naturally indexed by $NC(W)$, and
\item[$\bullet$] with an isomorphism of $\CC[W \times C]$-modules
\begin{equation}
\label{enlarged-NC(W)-isomorphism}
S/(\thetas) \underset{\CC[W \times C]-\text{mod}}{\cong} \CC[P]
\end{equation}
\end{enumerate}
where $W, C$ act on $S/(\thetas)$ by the usual linear, scalar substitutions
as in \eqref{Springer-isomorphism}?
\end{question}

A positive answer to this question would lead to a more conceptual proof of 
Theorem~\ref{main-result} along the lines of the proof of Theorem~\ref{Springer-CSP}:
restrict \eqref{enlarged-NC(W)-isomorphism} to the $W$-fixed subspaces as $\CC[C]$-modules, and 
then take the trace of any  element of $C$ to give the desired CSP.

\subsection{Speculation on the finite torus $Q/pQ$}
\label{finite-torus-section}
Here we restrict to the case where $W$ is a Weyl group, that is, a crystallographic
real reflection group.  In this case there is a root system for $W$ and the {\it root lattice} $Q$
given by the $\ZZ$-span of all roots is stabilized by $W$.
Haiman \cite[\S7]{Haiman} suggests considering, for various positive integers $p$,
the $W$-action on the ``finite torus'' $Q/pQ$, whose cardinality is $p^n$.  

\begin{proposition}\cite[Proposition 7.4.1]{Haiman}, \cite[Proposition 1.7]{BerestEtingofGinzburg},
\cite[\S 5 Theorem part 4]{Gordon}
Let $W$ be a Weyl group $W$ with root lattice $Q$, and $p$ any positive integer.
Then as a $W$-representation, the finite torus $Q/pQ$ has the same
character formula $\chi_{Q/pQ}(w) = p^{\dim_\CC(V^w)}$ as that
appearing in \eqref{V*-formulae}.

Consequently, when $\thetas$ is an hsop carrying $V(-p)(=V^*(-p))$, one has an isomorphism
of $\CC[W]$-modules
$$
S/(\thetas) \underset{\CC[W]-\text{mod}}{\cong} \CC[Q/pQ]
$$
and
$$
\dim_\CC ( S/(\thetas) )^W = \prod_{i=1}^n \frac{p+d_i-1}{d_i}.
$$

In particular, when $p=h+1$ and $\thetas$ is the hsop carrying $V(-p)$ described in
Remark~\ref{known-hsops}, one has
\begin{equation}
\label{parking-function-representation}
S/(\thetas)  \underset{\CC[W]-\text{mod}}{\cong} \CC[Q/(h+1)Q]
\end{equation}
and the dimension formula 
$$
\begin{aligned}
\dim_\CC ( S/(\thetas) )^W &= \dim_\CC \CC[Q/(h+1)Q]^W\\
                          &= | W\backslash Q/(h+1)Q |\\ 
                          &= \Cat(W).
\end{aligned}
$$
\end{proposition}

When $p=h+1$, these $W$-orbits $W \backslash Q/(h+1)Q$ in the
finite torus $Q/(h+1)Q$ counted in the last formula are
known (by work of Athanasiadis \cite{Athanasiadis} and of
Cellini, Papi, et al; see Armstrong \cite[\S 1.1]{Armstrong1})
to biject with various other objects of cardinality $\Cat(W)$, such as 
antichains in the poset of positive roots for $W$, or dominant chambers in the Shi hyperplane
arrangement for $W$.  

Our second speculation is an attempt to refine the $\CC[W]$-module isomorphism
in \eqref{parking-function-representation}, in the same way that
\eqref{Springer-isomorphism} refines \eqref{Chevalley-isomorphism}.

\begin{question}
\label{finite-torus-question}
Let $W$ be a Weyl group with root lattice $Q$, and let
$p=h+1$.  Can one naturally define a $C$-action on $\CC[Q/pQ]$ and choose the 
hsop $\thetas$ carrying $V(-p)$ in such a way that
\begin{enumerate}
\item[$\bullet$] 
the $C$-action on $\CC[Q/pQ]$ commutes with the $W$-action, leading to a
$W \times C$-action, and
\item[$\bullet$]
one has an isomorphism of $\CC[W \times C]$-modules
$$
S/(\thetas) \underset{\CC[W \times C]-\text{mod}}{\cong} \CC[Q/(h+1)Q]
$$
\end{enumerate}
where $W, C$ act on $S/(\thetas)$ by the usual linear, scalar substitutions
as in \eqref{Springer-isomorphism}?
\end{question}

If both Questions~\ref{enlarge-NC(W)-question} and Questions~\ref{finite-torus-question}
have affirmative answers, then the string of isomorphisms of $\CC[W \times C]$-modules
$$
\CC[P] 
\underset{\CC[W \times C]-\text{mod}}{\cong} 
S/(\thetas) 
\underset{\CC[W \times C]-\text{mod}}{\cong} 
\CC[Q/(h+1)Q],
$$
and the ensuing $\CC[C]$-module isomorphisms between their $W$-fixed subspaces of
dimension $\Cat(W)$, would provide a sought-after missing link between $NC(W)$
and $W\backslash Q/(h+1)Q$.

\section{Remarks and further questions}
\label{variations-section}

\subsection{$W$-$q$-Narayana and $W$-$q$-Kreweras numbers?}

Athanasiadis \cite{Athanasiadis} studied the refinement of the
$W$-Catalan numbers $\Cat(W)$ by the {\it $W$-Narayana numbers $\Nar(W,k)$};
these turn out to count (among other things) 
the elements $w$ of $NC(W)$ having fixed subspace $V^w$ of codimension $k$,
that is, the elements of rank $k$ in the poset structure on $NC(W)$.
In type $A_{n-1}$, there is a well-known  {\it $q$-Narayana number} $\Nar(n,k,q)$
that refines the $q$-Catalan number $\Cat(W,q)$, and 
\cite[Theorem 7.2]{RSW} actually proves a $CSP$ for these $q$-Narayana numbers
that refines Theorem~\ref{RSW-result}.  
Note that when the Coxeter element $c$ in $W$ acts by conjugation on $W$, it preserves
ranks (codimension of fixed spaces), and hence the cyclic group $C$ generated by $c$ acts
on each rank of $NC(W)$.

\begin{question}
Can one define $q$-Narayana numbers  $\Nar(W,k,q)$ for well-generated $W$
so that they exhibit the $CSP$ for this $C$-action on rank $k$ of $NC(W)$?  
\end{question}

One can be even more refined in counting elements of $NC(W)$.  In his original
study of the (type $A_{n-1}$) noncrossing partitions $NC(n)$, Kreweras \cite{Kreweras} showed
that the number of partitions in $NC(n)$ with exactly $m_i$ blocks of size $i$ 
for $i=1,2,\ldots$ and $k:=\sum_i m_i$ blocks total is 
\begin{equation}
\label{Kreweras-formula}
\frac{n(n-1)\cdots(n-k+2)}{m_1! m_2! \cdots}.
\end{equation}
We sketch below a (brute force) verification of the following more refined CSP.
\begin{theorem}
\label{noncrossings-by-type-CSP}
The triple $(X,X(q),C)$ given by
$$
\begin{aligned}
X&:=\text{ noncrossing partitions of }[n]\text{ with }m_i\text{ blocks of size }i\\
\label{q-Kreweras}
X(q)&:=\frac{[n]_q [n-1]_q \cdots [n-k+2]_q}{[m_1]!_q [m_2]!_q \cdots}\\
C&:=\ZZ_n \text{ via rotation}.
\end{aligned}
$$
exhibits the CSP.
\end{theorem}
\begin{proof}(Sketch)
For $d \geq 2$ dividing $n$ one can count the $d$-fold symmetric noncrossing partitions 
of $[n]$ with $m_i$ blocks of size $i$, as follows.  There can be at most one block which is 
itself $d$-fold symmetric-- two such blocks would cross each other.
Hence almost all of the multiplicities $m_i$ must be divisible by $d$, with at
most one such multiplicity $m_{j}$ congruent to $1$ modulo $d$, and furthemore this $m_j$ (if it exists)
must have its block size $j$ divisible by $d$.  
Define 
$$
n':=\frac{n}{d},\quad
j':=\frac{j}{d} \quad
k':=\left\lfloor \frac{k}{d} \right\rfloor,\quad
m'_i:=\left\lfloor \frac{m_i}{d} \right\rfloor.
$$
Relabel the numbers $[2n']=\{1,2,\ldots,2n\}$ as $1,2,\ldots,n,-1,-2,\ldots,-n$, and then
restrict the blocks of the $d$-fold symmetric partition of $[n]$ to this set $[2n']$;
if there is a (unique) $d$-fold symmetric $j$-block, then call its restriction the ``zero block''.
This gives a bijection to type $B_{n'}$ noncrossing partitions considered in \cite{Reiner} having
$m'_i$ nonzero blocks of size $i$, of which there are
\begin{equation}
\label{Athaniasiadis-formula}
\frac{n'(n'-1)\cdots(n'-k'+1)}{m'_1! m'_2! \cdots}
\end{equation}
by a formula of 
Athanasiadis \cite[Theorem 2.3]{Athanasiadis-block-sizes}.

We must check this agrees with the evaluation $X(\omega)$ where $\omega$ is a primitive
$d^{th}$ root-of-unity.  Since $d$ divides $n$, the multiplicity of $\omega$ as a root
in the numerator of $X(q)$ is exactly 
$
\left\lceil \frac{k-1}{d} \right\rceil \geq \frac{k}{d}-\frac{1}{d},
$ 
while in the denominator its multiplicity is exactly 
$
\sum_i m'_i = \frac{k}{d} - \sum_i \frac{m''_i}{d}
$
where $m''_i$ is the (nonnegative) remainder of $m_i$ on division by $d$.
Hence $X(\omega)$ will vanish unless at most one of the $m''_i$ is nonzero,
and that unique nonzero $m''_j=1$.  In other words, $X(\omega)$ vanishes unless
the $m_i$ satisfy exactly the same congruence conditions as required to have
a $d$-fold symmetric noncrossing partition.
Under those conditions, one can evaluate $X(\omega)$ by
matching up numerator factors and
denominator factors as $\frac{[a]_q}{[b]_q}$ with $a \equiv b \mod d$, and
use Proposition~\ref{numerical-prop}; one finds that $X(\omega)$ agrees with
\eqref{Athaniasiadis-formula}.
\end{proof}

One can generalize Kreweras' formula \eqref{Kreweras-formula} from type $A_{n-1}$ to all well-generated
complex reflection groups:  one counts elements $w$ in $NC(W)=[e,c]$ from a given $W$-conjugacy class.
This turns out to be equivalent to counting the elements $w$ in $NC(W)$ whose fixed subspace $V^w$ lies within
a given $W$-orbit $W \cdot V'$ of subspaces $V'$.  
When $W$ happens to be a Weyl group, it has been verified in a case-by-case fashion 
\cite[Theorem 6.3]{AthanasiadisReiner} that the above elements
are equinumerous with the  $W$-orbits in $Q/(h+1)Q$ for which the $W$-stabilizer of any member of the class
is $W$-conjugate to the (pointwise) $W$-stabilizer of $V'$.  Equivalently, this is the
number of order filters in the poset of positive roots for which the
minimal roots of the filter have common perpendicular space lying in the $W$-orbit $W \cdot V'$.

When $W$ is Weyl group, there is a generalization of the Kreweras formula \eqref{Kreweras-formula} that
counts these sets, due to Sommers and Trapa, Douglass, and Broer (see \cite[\S 6]{Sommers}) 
They consider the {\it Orlik-Solomon exponents} $e_1,\ldots,e_\ell$ for the hyperplane arrangement
one gets by restricting the hyperplanes of $W$ to $V'$, and show that there are exactly
\begin{equation}
\label{Orlik-Solomon-formula}
\frac{1}{[N_W(W'):W']} \prod_{i=1}^\ell (h+1-e_i)
\end{equation}
elements in any of the above sets. Here $W'$ is the pointwise stabilizer any element
in the $W$-orbit, while $N_W(W')$ is the normalizer of $W'$.

\begin{question}
Let $W$ be a Weyl group acting on $V$, and $C=\langle c \rangle$ for some Coxeter element
used to define $NC(W)$.  Let $V'$ be a subspace of $V$, with pointwise stabilizer $W'$.  
Let $X$ be a subset of elements $w$ in $NC(W)$ whose fixed subspace $V^w$ lies
in the $W$-orbit $W \cdot V'$.

Can one produce a $q$-analogue $X(q)$ of the formula \eqref{Orlik-Solomon-formula} so that the triple
$(X,X(q),C)$ exhibits the CSP?
\end{question}

\noindent
This can easily be done in types $A$ and $B$, and there is hope that the work
of Broer \cite{Broer} can shed some light on how to proceed uniformly.

\subsection{$W$-$q$-Fuss-Catalan numbers}

The $W$-$q$-Fuss-Catalan number
$$
\Cat^m(W,q):= \prod_{i=1}^{n} \frac{ [m h+d_i]_q }{ [d_i]_q }.
$$
that appeared in Theorem~\ref{Cat-consequences} relates to 
work of Armstrong \cite[Chapter 3]{Armstrong1}.  
He has defined the {\it $m$-divisible} noncrossing
partitions $NC^m(W)$ for any real reflection group $W$, having
cardinality $|NC^m(W)| = \Cat^m(W,1)$:  they are the $(m+1)$-tuples
$(w_0,w_1,\ldots,w_m)$ factoring $c=w_0 w_1 \cdots w_m$ in such a way that the
lengths are additive:  $\sum_{i=0}^m \ell_R(w_i) = \ell_R(c) ( = |S|)$.
In the case $m=1$, the map $(w_0,w_1) \mapsto w_1$ identifies $NC^m(W)$ with $NC(W)$.
The set $NC^m(W)$ also has an extremely well-behaved poset structure:  it turns out
that for any element $(w_0,w_1,\ldots,w_m)$ in $NC^m(W)$, each $w_i$ lies in $NC(W)$, 
and one orders them by the product partial ordering $NC(W)^m$ on the last $m$-factors 
$(w_1,\ldots,w_m)$, forgetting $w_0$.

Armstrong \cite[Definition 3.4.15]{Armstrong1} also notes that
there is a natural cyclic group $C$ of order $mh$
acting on $NC^m(W)$ (generalizing the case $m=1$):  the generator of $C$ acts by
\begin{equation}
\label{armstrong-action}
\begin{matrix}
         & (w_0, & w_1,       & w_2, & w_3, & \ldots, & w_m)\\
 \mapsto & (v,   & cw_mc^{-1},& w_1, & w_2, & \ldots, &w_{m-1})
\end{matrix}
\end{equation}
where $v:=(c w_m c^{-1})w_0(c w_m c^{-1})^{-1}$.

Note that these definitions all make sense more generally for any well-generated complex
reflection group $W$.

\begin{conjecture}(\cite[Conjecture 5.4.7]{Armstrong1})
\label{conjecture-armstrong}
For a well-generated complex reflection group $W$ and any $m \geq 1$,
the triple $(X,X(q),C)$ given by
$$
\begin{aligned}
X &= NC^m(W) \\
X(q) &= \Cat^m(W,q) \\
C&=\ZZ/mh\ZZ
\end{aligned}
$$
exhibits the CSP.
\end{conjecture}

One can also combine this conjecture with the previous question, in that there are
well-defined {\it $m$-extended $W$-Narayana numbers} $\Nar^m(W,k)$, which count
the elements in $NC^m(W)$ of rank $k$ (and have other interpretations;  see 
Athanasiadis \cite{Athanasiadis}, Fomin and
Reading \cite[\S 7]{FominReading} and Armstrong \cite[Chapter 3]{Armstrong1}).  Is there a $q$-version
$\Nar^m(W,k,q)$ of these numbers, and an associated CSP for the cyclic group $C=\ZZ/mh\ZZ$
described above?

There is a variant of Conjecture \ref{conjecture-armstrong} involving
a modified action. Here, the generator of the cyclic
group $C'$ of order $(m+1)h$ acts on 
$NC^m(W)$ by
\begin{equation}
\label{bessis-action}
\begin{matrix}
         & (w_0, & w_1,       & w_2, & w_3, & \ldots, & w_m)\\
 \mapsto & (cw_mc^{-1},   & w_0,& w_1, & w_2, & \ldots, &w_{m-1})
\end{matrix}
\end{equation}

\noindent
One can readily check that if one includes $NC^m(W) \hookrightarrow NC^{m+1}(W)$ as
follows
$$
(w_0, w_1, w_2,, \ldots, w_m)
\mapsto 
(e, w_0, w_1, w_2, \ldots, w_m)
$$
then the cyclic action on $NC^{m+1}(W)$ defined by \eqref{armstrong-action}
restricts to the subset $NC^m(W)$, and coincides with the cyclic action defined in \eqref{bessis-action}.
The latter cyclic action on $NC^m(W)$ may also be interpreted in terms of the \emph{Garside automorphism}
of the \emph{$(m+1)$-divided category} associated with the dual braid
monoid of $W$ (see \cite{Bessis3}).

\begin{conjecture}
\label{conjecture-bessis}
For a well-generated complex reflection group $W$ and any $m \geq 0$,
the triple $(X,X(q),C')$ given by
$$
\begin{aligned}
X &= NC^m(W) \\
X(q) &= \Cat^m(W,q)\\
C'&=\ZZ/(m+1)h\ZZ
\end{aligned}
$$
exhibits the CSP.
\end{conjecture}

\noindent
Even the special case where $m=1$ in Conjecture~\ref{conjecture-bessis} is interesting, and
an open question, which strictly generalizes Theorem~\ref{main-result}.  
Here the generator of $C'=\ZZ/2h\ZZ$ acting on $NC^1(W)=NC(W)$
sends $w \mapsto cw^{-1}$, a map sometimes called the {\it Kreweras complementation},
and whose square is the generator $w \mapsto cwc^{-1}$ of the cyclic group $C=\ZZ/h\ZZ$
in Theorem~\ref{main-result}.

\begin{example} \rm \ 
When $W=\symm_n$ is of type $A_{n-1}$, so that $h=n$,
Conjecture~\ref{conjecture-bessis} has been verified for $m=1$ previously by 
D. White \cite{White-personalcommunication} via direct calculations.
In this case for arbitrary $m$, the cyclic action on $X=NC^m(W)$ is 
equivalent to the following:  $X$ may be identified
with the subset of the noncrossing
partitions $NC((m+1)n)$ which decompose $[(m+1)n]=\{1,2,\ldots,(m+1)n\}$ into 
exactly $n$ blocks each of size $m+1$, and
the cyclic group $C'=\ZZ/(m+1)n\ZZ$ acts on such partitions 
by cycling the numbers in $[(m+1)n]$ modulo $(m+1)n$.
\end{example}

\begin{example}\rm \ 
If $W$ is a Weyl group of type $E_8$, having degrees 
$$
2,8,12,14,18,20,24,30,
$$
$NC^1(W)$ admits an action of the cyclic group $C'$ of order $60$;
if $\phi\in C'$ is an element of order $4$, one checks that
$|NC^1(W)^{\phi}|= 88 = \Cat^1(W,\zeta_4)$. This number may be interpreted
%
%
as the number of objects in a \emph{dual Garside category} associated
with the (badly-generated) complex reflection group $G_{31}$, the centralizer
of a $4$-regular element in $W$ (\cite{Bessis2,Bessis3}). 
\end{example}

\subsection{Panyushev's cyclic action}

  Let $W$ be a crystallographic Weyl group for a root system $\Phi_W$, and
make a choice of positive roots $\Phi_W^+$.  There is a usual partial ordering on
$\Phi_W^+$ in which one $\alpha \leq \beta$ for two roots $\alpha, \beta$ if the
difference $\beta-\alpha$ is a nonnegative combination of positive roots.
It was mentioned in Section~\ref{finite-torus-section} that the collection $\antichains(\Phi_W^+)$ of {\it antichains} 
(= subsets containing only pairwise incomparable elements) for this ordering on $\Phi_W^+$
is a $W$-Catalan object in the sense that $|\antichains(\Phi_W^+)|=\Cat(W)$; see \cite[Chapter 1]{Armstrong1}.

   In this context, Panyushev \cite{Panyushev} recently defined the following cyclic action on 
$\antichains(\Phi_W^+)$:  the generator sends an antichain $A$ in $\Phi^+$ to the set of
maximal elements of $\Phi^+$ which lie above no
element of $A$.  He then makes several interesting conjectures about the orbit structure for
this cyclic action.  In particular he conjectures that the order of the cyclic group it generates
always divides $2h$ where $h$ is the Coxeter number.  We add to his conjectures the following
conjecture of our own.

\begin{conjecture}
\label{Panyushev-action-gives-CSP}
For any crystallographic Weyl group $W$, with positive roots $\Phi^+$, the triple $(X,X(q),C)$ given by
$$
\begin{aligned}
X &= \antichains(W) \\
X(q) &= \Cat(W,q)\\
C&=\ZZ/2h\ZZ, \text{ generated by Panyushev's action.}
\end{aligned}
$$
exhibits the CSP.
\end{conjecture}

\noindent
If one assumes the $m=1$ case of Conjecture~\ref{conjecture-bessis}, 
then Conjecture~\ref{Panyushev-action-gives-CSP} is
equivalent to the conjecture that there exists a $\ZZ/2h\ZZ$-equivariant bijection 
$$
\antichains(W) \rightarrow NC(W)
$$
in which the action on $\antichains(W)$ is Panyushev's, and the action on $NC(W)$
is the Kreweras complementation.  Presumably there also exists an ``$m$-version'' 
of Conjecture~\ref{Panyushev-action-gives-CSP},
involving an extension of Panyushev's action to the $m$-version
of antichains counted by $\Cat^m(W,q)$ which are described by Athanasiadis in \cite[Theorem 1.2]{Athanasiadis},
and involving $X(q)= \Cat^m(W,q)$.

\begin{remark}
The CSP Conjectures~\ref{conjecture-armstrong} and \ref{conjecture-bessis} 
both involve $X(q)=\Cat^m(W,q)$ and $X=NC^m(W)$ (but with two different cyclic groups, of orders
$mh$ and $(m+1)h$, respectively).  

On the other hand, the discussion surrounding Conjecture~\ref{Panyushev-action-gives-CSP}
involves $X(q)=\Cat^m(W,q)$ but with $X$ given by the $m$-version of antichains of $\Phi^+_W$.

We last mention another CSP involving $X(q)=\Cat^m(q)$ and with $X$ given by
the third known class of $W$-Catalan objects: the set of {\it clusters} in the 
generalized cluster complexes of Fomin and Reading \cite{FominReading}.
Here the cyclic group $C$ is of order $mh+2$, generated by the action on clusters 
of the so-called {\it deformed Coxeter element}.  A CSP here was conjectured by the
authors of \cite{RSW}, appeared as \cite[Conjecture 5.4.8]{Armstrong1},
and was recently verified directly by Eu and Fu \cite[Theorem 1.2 and Conjecture 8.1]{EuFu}.
\end{remark}

\section*{Acknowledgements}
The authors thank Drew Armstrong, Christos Athanasiadis, 
Sen-Peng Eu, Dennis Stanton, Dennis White, and two anonymous referees for 
helpful comments.  
The second author particularly thanks Stephen Griffeth for helpful tutorials 
on rational Cherednik algebras. 

Both authors thank the American Institute of Mathematics for their 
support and hospitality during the workshop on
``Braid groups, clusters, and free probability'' in January 2005.

\end{document}